\documentclass[aps,prl,reprint,groupedaddress]{revtex4-2}
\usepackage{amsmath}
\usepackage{bm}
\usepackage{amsfonts}
\usepackage{latexsym}
\usepackage{amssymb}
\usepackage{stmaryrd}
\usepackage{overpic}
\usepackage{graphicx}
\usepackage{multirow}
\usepackage{mathrsfs}
\usepackage{float}
\usepackage{color}
\usepackage{setspace}
\usepackage{subfigure}
\usepackage{dsfont}

\begin{document}

\title{A Novel Computational Method for Band Structures of Dispersive Photonic Crystals}
\author{Wenqiang Xiao$^1$ and Jiguang Sun$^2$}
\affiliation{$^1$Beijing Computational Science Research Center, Beijing 100193, China.\\
$^2$Michigan Technological University, Houghton, MI 49931, U.S.A.}

\begin{abstract}
We propose a new method to compute band structures of dispersive photonic crystals. It can treat arbitrarily frequency-dependent, lossy or lossless materials. The band structure problem is first formulated as the eigenvalue problem of an operator function.  Finite elements are then used for discretization. Finally, the spectral indicator method is employed to compute the eigenvalues. Numerical examples in both the TE and TM cases are presented to show the effectiveness. There exist very few examples in literature for the TM case and three examples in this paper can serve as benchmarks.
\end{abstract}

\maketitle
\section{Introduction}
Photonic crystals are periodic structures composed of dielectric materials.
Due to the presence of band gaps, they have the capability to control and manipulate electromagnetic wave propagation
with numerous applications in optical communications, filters, lasers, switches and optical transistors
\cite{Soukoulis, Sakoda2001,YangEtal}.
Effective calculations of the band
structures are important to identify novel optical phenomena and develop new devices.

If the permittivity is independent of the frequency, the band structure calculation needs to solve linear eigenvalue problems of the Maxwell's equations. 
There exist many successful methods such as 
the plane wave expansion \cite{ho1990}, the finite-difference time-domain (FDTD) method \cite{qiu2000}, 
and the finite element method \cite{axmann1999, dobson1999}, etc. 
In contrast, frequency dependent permittivities usually lead to nonlinear eigenvalue problems. 
Much less effective methods have been developed, e.g., the
generalizations of the plane wave expansion method, the FDTD method, 
the cutting plane method, 
and path-following algorithm \cite{Soukoulis, spence2005, kuzmiak1998,ito2001,toader2004}.
These methods either solve the nonlinear eigenvalue problems using the Newton's iteration, which requires accurate initial values, 
or impose certain restrictions on the material property, e.g., the electric permittivity has some special forms or the shape of the crystal
needs to be regular. Nonetheless, a general method for effective band structure calculation for dispersive photonic crystals is highly desirable.

In this study, we propose a novel approach to compute the band structures of photonic crystals of dispersive materials.
The permittivity can depend arbitrarily on the frequency as well as the position. The materials can be lossy or lossless. 
The approach works for both TE case and TM case and can treat crystals with irregular shapes. 
Roughly speaking, the band structure problem is formulated as the eigenvalue problem of an operator-valued function \cite{karma1996II, gohberg2009, engstrom2010, xiao2020}. 
Then finite elements are used to discretize the operator-valued function.
Finally, a spectral indicator method (SIM) is developed to practically compute the eigenvalues of the operator-valued function.
This version of SIM extends the idea in \cite{huang2016, huang2018, Huang2020, SunZhou2016, xiao2020} for the generalized eigenvalue problems of non-Hermitian matrices and is particularly suitable to compute eigenvalues of operator-valued functions. 

\section{Eigenvalues of an Operator-Valued Function}
Photonic band structures are determined by the eigenvalue problems of the Maxwell's equations \cite{axmann1999, raman2010}
\begin{equation}\label{maxwell}
\left\{
\begin{aligned}
& \nabla\times\frac{1}{\epsilon(x,\omega)}\nabla\times\mathbf{H}=\left(\frac{\omega}{c}\right)^2\mathbf{H},\\
& \nabla\cdot\mathbf{H}=0,
\end{aligned}
\right.
\end{equation}
where $\mathbf{H}$ is the magnetic field, $\omega$ is the frequency, $\epsilon(x,\omega)$ is the electric permittivity, and $c$ is the speed of light in the vacuum.
Let $\Omega$ be a compact set on the complex plane $\mathbb{C}$. 
In 2D, for $\omega \in \Omega$, the Maxwell's equations \eqref{maxwell} can be reduced to the transverse electric (TE) case 
\begin{equation}\label{TE}
-\Delta\psi=\left(\frac{\omega}{c}\right)^2\epsilon(x,\omega)\psi.
\end{equation}
or the transverse magnetic (TM) case:
\begin{equation}\label{TM}
-\nabla\cdot\left(\frac{1}{\epsilon(x,\omega)}\nabla \psi\right)=\left(\frac{\omega}{c}\right)^2\psi.
\end{equation}

\begin{figure}[h!]
\begin{center}
{ \scalebox{0.40} {\includegraphics{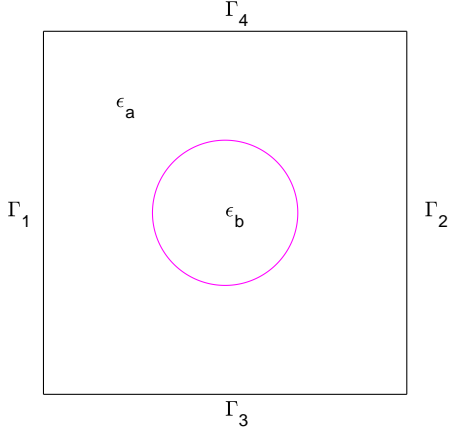}}}
\caption{A square cell $D_0$ with a disc of radius $r$ at its center.}
 \label{geometry}
\end{center}
\end{figure}

For simplicity, we assume that the photonic crystal has the unit periodicity on a square lattice (Fig.~\ref{geometry}). 
Since finite element methods will be used for the discretization of the equation, other lattices/shapes can be treated in the same way.
Define
$\mathbb{Z} =\{0,\pm 1,\pm 2,\cdots\}$ and the lattice $\Lambda=\mathbb{Z}^2$. The electric permittivity is a periodic function given by
\begin{equation*}
\epsilon(x+n,\omega)=\epsilon(x,\omega), ~~x=(x_1, x_2)\in\mathbb{R}^2, n\in\Lambda.
\end{equation*}
The infinite periodic domain $D:=\mathbb{R}^2/\mathbb{Z}^2$
can be identified with the unit square $D_0:=(0, 1)^2$ by imposing periodic boundary conditions.
Introduce the quasimomentum vector ${\boldsymbol k}\in \mathcal{K}$, where
$$\mathcal{K}=\{{\boldsymbol k}\in \mathbb{R}^2|-\pi\leq k_j\leq\pi,j=1,2\}$$
is the Brillouin zone.
Using the Floquet transform \cite{kuchment1993}, \eqref{TE} and \eqref{TM}, respectively,
can be written as eigenvalue problems parametrized by ${\boldsymbol k}$:
\begin{equation}\label{TE-equation}
 -(\nabla+i {\boldsymbol k})\cdot(\nabla+i{\boldsymbol k})u(x)=\left(\frac{\omega}{c}\right)^2\epsilon(x,\omega)u(x)~~ {\rm in}~D_0
 \end{equation}
 and
 \begin{equation}\label{TM-equation}
 -(\nabla+i{\boldsymbol k})\cdot\frac{1}{\epsilon(x,\omega)}(\nabla+i{\boldsymbol k})u(x)=\left(\frac{\omega}{c}\right)^2 u(x)~~{\rm in}~D_0.
\end{equation}
Here $u$ is the Floquet transform of $\psi$, a periodic function in both $x_1$ and $x_2$.
Let $\Gamma_i , i=1,\cdots,4,$ be the four parts of $\partial D_0$ (Fig.~\ref{geometry}).
\begin{equation}\label{BCs}
\left\{
\begin{aligned}
&u|_{\Gamma_1}=u|_{\Gamma_2},~~u_{x_1}|_{\Gamma_1}=u_{x_1}|_{\Gamma_2},\\
&u|_{\Gamma_3}=u|_{\Gamma_4},~~u_{x_2}|_{\Gamma_3}=u_{x_2}|_{\Gamma_4}.
\end{aligned}
\right.
\end{equation}

For any $\omega\in \Omega$, assume that the permittivity $\epsilon(x,\omega)$ is a piecewise continuous function satisfying
\begin{equation}\label{boundness}
0<C_0\le |\epsilon(x,\omega)|\le C_1<\infty,~x\in D_0.
\end{equation}
 
Let $L_2(D_0)$ be the space of square-integrable functions and $H^1(D_0)$ be the space of functions in $L^2(D)$ whose first order partial derivatives
are also square-integrable.
Define the subspace of functions in $H^1(D_0)$ with periodic boundary conditions by
\begin{equation}
H_p^1(D_0)=\{v\in H^1(D_0)|~v|_{\Gamma_1}=v|_{\Gamma_2},v|_{\Gamma_3}=v|_{\Gamma_4}\}.
\end{equation}
Multiplying \eqref{TE-equation} and \eqref{TM-equation} by a test function $v\in H_p^1(D_0)$ and integrating by parts, the weak formulations are to find 
$(\omega, u) \in \Omega \times H_p^1(D_0)$ such that
\begin{equation}\label{TE weak form}
\int_{D_0}(\nabla+i{\boldsymbol k})u\cdot\overline{(\nabla+i{\boldsymbol k})v}dx=\left(\frac{\omega}{c}\right)^2\int_{D_0} \epsilon(x,\omega)u\overline{v}dx
\end{equation}
for all  $v\in H_p^1(D_0)$ in the TE-case and
\begin{equation}\label{TM weak form}
\int_{D_0}\frac{1}{\epsilon(x,\omega)}(\nabla+i{\boldsymbol k})u\cdot\overline{(\nabla+i{\boldsymbol k})v}dx=\left(\frac{\omega}{c}\right)^2\int_{D_0} u\overline{v}dx
\end{equation}
for all  $v\in H_p^1(D_0)$ in the TM-case.
If $\epsilon(x, \omega)$ depends on $\omega$,  \eqref{TE weak form} and \eqref{TM weak form} are nonlinear eigenvalue problems in general.

Using \eqref{TE weak form} and \eqref{TM weak form}, we define an operator-valued function 
\[
T: \Omega \to \mathcal{L}(H_p^1(D_0), H_p^1(D_0)),
\] 
where $\mathcal{L}(H_p^1(D_0), H_p^1(D_0))$
denotes the space of bounded linear operators from $H_p^1(D_0)$ to $H_p^1(D_0)$.
For the TE case, $T(\omega):=T_{TE}(\omega)$ such that
\begin{eqnarray}\nonumber
&& \left(T_{TE}(\omega)u,v\right) \\ \nonumber
&=& \int_{D_0} \nabla u \cdot \nabla \overline{v}+i{\boldsymbol k} u \cdot \nabla \overline{v} -i\nabla u \cdot {\boldsymbol k}\overline{v}+|{\boldsymbol k}|^2u\overline{v} \\
\label{bTE} && \qquad -  \left(\frac{\omega}{c}\right)^2\epsilon(x, \omega)u\overline{v}dx  ~\text{for all } v \in H_p^1(D_0).
\end{eqnarray}
Similarly, for the TM case, $T(\omega):=T_{TM}(\omega)$ such that
\begin{eqnarray}\nonumber
&& \left(T_{TM}(\omega)u,v\right) \\
\nonumber &=&\int_{D_0} \frac{1}{\epsilon(x, \omega)} \left(\nabla u \cdot \nabla \overline{v}+i{\boldsymbol k} u \cdot \nabla \overline{v} -i\nabla u \cdot {\boldsymbol k}\overline{v}+|{\boldsymbol k}|^2u\overline{v}\right) \\
\label{TM eqn} && \qquad -  \left(\frac{\omega}{c}\right)^2u\overline{v}dx ~\text{for all } v \in H_p^1(D_0).
\end{eqnarray}
Then the eigenvalue problems \eqref{TE weak form} and \eqref{TM weak form} are equivalent to the operator eigenvalue problem of finding
$\omega \in \Omega$ and $u \in H_p^1(D_0)$ such that
\begin{equation}\label{Tomegau}
T(\omega)u = 0.
\end{equation}

\section{FEM Discretization and Spectral Indicator Method}
In this section, we employ the finite element mehtods to discretize $T(\omega), \omega \in \Omega$,
and propose a spectral indicator method to compute the eigenvalues of \eqref{Tomegau}. 

Let $\mathcal{T}_h$ be a regular triangular mesh for $D_0$, where $h$ is the mesh size.
Let $V_h\subset H_p^1(D_0)$ be the associated Lagrange finite element space
and $\phi_i, i=1,\cdots,N$, be the basis functions for $V_h$.
For a fixed ${\boldsymbol k} \in\mathcal{K}$, let 
\[
M^h_\epsilon,~~ M^h, ~~S^h, ~~A^{1,h}, ~~ A^{2,h}
\]
be the matrices corresponding to the terms 
\begin{eqnarray*}
&&(\epsilon(x, \omega) u_h, \overline{v}_h), \\ 
&&(u_h, \overline{v}_h), ~(\nabla u_h, \nabla \overline{v}_h), \\
&&({\boldsymbol k}u_h, \nabla \overline{v}_h),\\ 
&&(\nabla u_h, {\boldsymbol k}\overline{v}_h)
\end{eqnarray*}
in \eqref{bTE}, 
and 
\[
M_{\epsilon^{-1}}^h, ~~S_{\epsilon^{-1}}^h, ~~A_{\epsilon^{-1}}^{1,h}, ~~A_{\epsilon^{-1}}^{2,h}
\]
be the matrices corresponding to the terms
\begin{eqnarray*}
&&(\epsilon(x,\omega)^{-1} u_h, \overline{v}_h),\\
&&(\epsilon(x,\omega)^{-1}\nabla u_h, \nabla \overline{v}_h), \\
&&(\epsilon(x,\omega)^{-1}{\boldsymbol k}u_h, \nabla \overline{v}_h), \\
&&(\epsilon(x,\omega)^{-1}\nabla u_h, {\boldsymbol k}\overline{v}_h)
\end{eqnarray*} 
in \eqref{TM eqn}.

In matrix form, the finite element discretization of \eqref{Tomegau} can be written as the problem of finding $\omega$ and a vector $\vec{u}$ such that
\begin{equation}\label{Thou0}
{\mathbb T}^h(\omega) \vec{u} = \vec{0},
\end{equation}
where, for the TE-case, ${\mathbb T}^h(\omega):={\mathbb T}_{TE}^h(\omega)$ is given by
\[
{\mathbb T}_{TE}^h(\omega)=S^h+iA^{1,h}-iA^{2,h}+|{\boldsymbol k}|^2 M^h-\left(\frac{\omega}{c}\right)^2 M^h_\epsilon
\]
and, for the TM-case, ${\mathbb T}^h(\omega):={\mathbb T}_{TM}^h(\omega)$ is given by
\[
{\mathbb T}_{TM}^h(\omega)=
S_{\epsilon^{-1}}^h+iA_{\epsilon^{-1}}^{1,h}-iA_{\epsilon^{-1}}^{2,h}+|{\boldsymbol k}|^2 M_{\epsilon^{-1}}^h-\left(\frac{\omega}{c}\right)^2 M^h.
\]

Now we propose a spectral indicator method to compute all the eigenvalues of ${\mathbb T}^h$ in a compact region $\Omega$ of interests on $\mathbb C$.
Without loss of generality, let $\Omega\subset\mathbb{C}$ be a square and $\Theta$ be the circle circumscribing 
$\Omega$. For example, $\Omega$ can be a square which is symmetric with respect to the real axis if one wants to find eigenvalues close to the real axis.

When ${\mathbb T}^h(\omega)$ is holomorphic in $\Omega$ and ${\mathbb T}^h(\omega)^{-1}$ exists and is bounded for all $\omega \in\Theta$, 
one can define an operator $\Pi$ by
\begin{equation}\label{Pi}
\Pi=\frac{1}{2\pi i}\int_\Theta {\mathbb T}^h(\omega)^{-1} d\omega.
\end{equation}
Let $\vec{g}_h$ be a random vector. If ${\mathbb T}^h$ has no eigenvalues in $\Omega$, $\Pi \vec{g}_h=0$. 
On the other hand, if ${\mathbb T}^h$ has at least one eigenvalue in $\Omega$, $\Pi \vec{g}_h\neq 0$ almost surely. 
This is the key idea of the spectral indicator method (see \cite{huang2016, huang2018, Huang2020}).

Using the trapezoidal rule to 
approximate $\Pi \vec{g}_h$, we define an indicator for $\Omega$ by
\begin{equation}\label{Iomega}
I_\Omega := \left|\Pi \vec{g}_h \right| \approx \bigg|\frac{1}{2\pi i}\sum\limits_{j=1}^{m_0}w_j \vec{x}_h(\omega_j)\bigg|,
\end{equation}
where $m_0$ is the number of quadrature points, $\omega_j$ are the quadrature points and $w_j$ are the weights. 
Here $\vec{x}_h(\omega_j)$ is the solution of the linear system ${\mathbb T}^h(\omega_j) \vec{x}_h(\omega_j)=\vec{g}_h$,
avoiding the computation of ${\mathbb T}^h(\omega)^{-1}$ in \eqref{Pi}.

The indicator $I_\Omega$ is used to decide if $\Omega$ contains eigenvalues or not. 
Let $\delta_0>0$ be the threshold value. 
If $I_{\Omega}>\delta_0$, $\Omega$ is said to be admissible and divided into small squares.
Otherwise, $\Omega$ is not admissible and contains no eigenvalues.
One computes the indicators of the small squares and decide if they are admissible.
The procedure continues until the diameter of an admissible square $d(\Omega)$ is less than $\beta_0$. Consequently, eigenvalues are identified with precision $\beta_0$.

For a fixed ${\boldsymbol k} \in \mathcal{K}$, the following algorithm SIM-H (spectral indicator method for holomorphic operator-valued functions) 
computes all the eigenvalues of ${\mathbb T}^h$ in a collection of uniform square regions $\{\Omega_k^0, k=1, \ldots, I^0\}$ of interests.
We shall use the TE case as an example. The TM case is similar.

\begin{itemize}
\item[] {\bf SIM-H:}
\item[-] Generate a a triangular mesh $D_0$. 
\item[-] Given a collection of square regions $\{\Omega_k^0\}$.
\item[-] Choose  the precision $\beta_0$ and the threshold $\delta_0$.
\item[1.] Assemble the matrices  $M^h, S^h, A^{1,h}, A^{2,h}$.
\item[2.] Choose a random $\vec{g}_h$ and set $i=0$.
\item[3.] At level $i$, do
	\begin{itemize}
	\item For each square $\Omega^i_k$ at current level $i$, evaluate the indicator $I_{\Omega^i_k}$ as follows.
		\begin{itemize}
			\item At each quadrature point $\omega_j$, assemble $M^h_\epsilon(\omega_j)$ and ${\mathbb T}^h(\omega_j)$.
			\item Solve $\vec{x}_h(\omega_j)$ for ${\mathbb T}^h(\omega_j)  \vec{x}_h(\omega_j) = \vec{g}_h$.
			\item Compute
		\[
			I_{\Omega^i_k} :=  \bigg|\frac{1}{2\pi i}\sum\limits_{j=1}^{m_0}w_j \vec{x}_h(\omega_j)\bigg|.
		\]
		\end{itemize}
	\item If $ I_{\Omega^i_k} > \delta_0$, mark $\Omega^i_k$ as admissible. Otherwise, discard $\Omega^i_k$.
	\item If  $d(\Omega_k^i) > \beta_0$, uniformly divide all the admissible squares into smaller squares and leave them to the next level $i \leftarrow i+1$.
		Otherwise, stop.
	\end{itemize}
\item[4.] Output the eigenvalues (centers of the admissible square regions).
\end{itemize}
The algorithm is highly scalable. The collection of squares $\{\Omega_k^i, k=1, \ldots, I^i\}$ at any level $i$ can be treated parallelly.

\section{Numerical Examples}
We present several examples by showing the dispersion relations $\omega({\boldsymbol k})$ with ${\boldsymbol k}$ 
moving from $\Gamma = (0, 0)$,  to $X = (\pi, 0)$, to $M=(\pi, \pi)$, and back to $\Gamma$.
Consider an array of identical, infinitely long, parallel cylinders, embedded in vacuum, 
whose intersection with a perpendicular plane forms a simple square lattice with a disc inside.
In all the examples, $\epsilon_a=1$ is the permittivity of the vacuum and
$\epsilon_b$ is the permittivity in the disc.
Define the filling fraction $f=\pi r^2/a^2$, where $a=1$ is the side length of the unit cell $D_0$, $r$ is the radius of the disc at the center of $D_0$. 
The linear Lagrange element is used for discretization.
In all numerical examples, we set $\delta_0=0.01$ and $\beta_0=10^{-4}$. We note that $\delta_0$ is usually problem dependent and can be chosen by test and error.

When materials are lossy, the eigenvalues are complex in general. 
In the following examples, the imaginary parts are relative small and, the real parts of the eigenvalues are used in the dispersion diagrams.
There exist very few examples in literature for the TM case, which are deemed to be challenging. We present one dispersion diagram for the TE case and three for the TM case. 

\begin{figure}[h!]
\begin{center}
\begin{tabular}{lll}
\resizebox{0.42\textwidth}{!}{\includegraphics{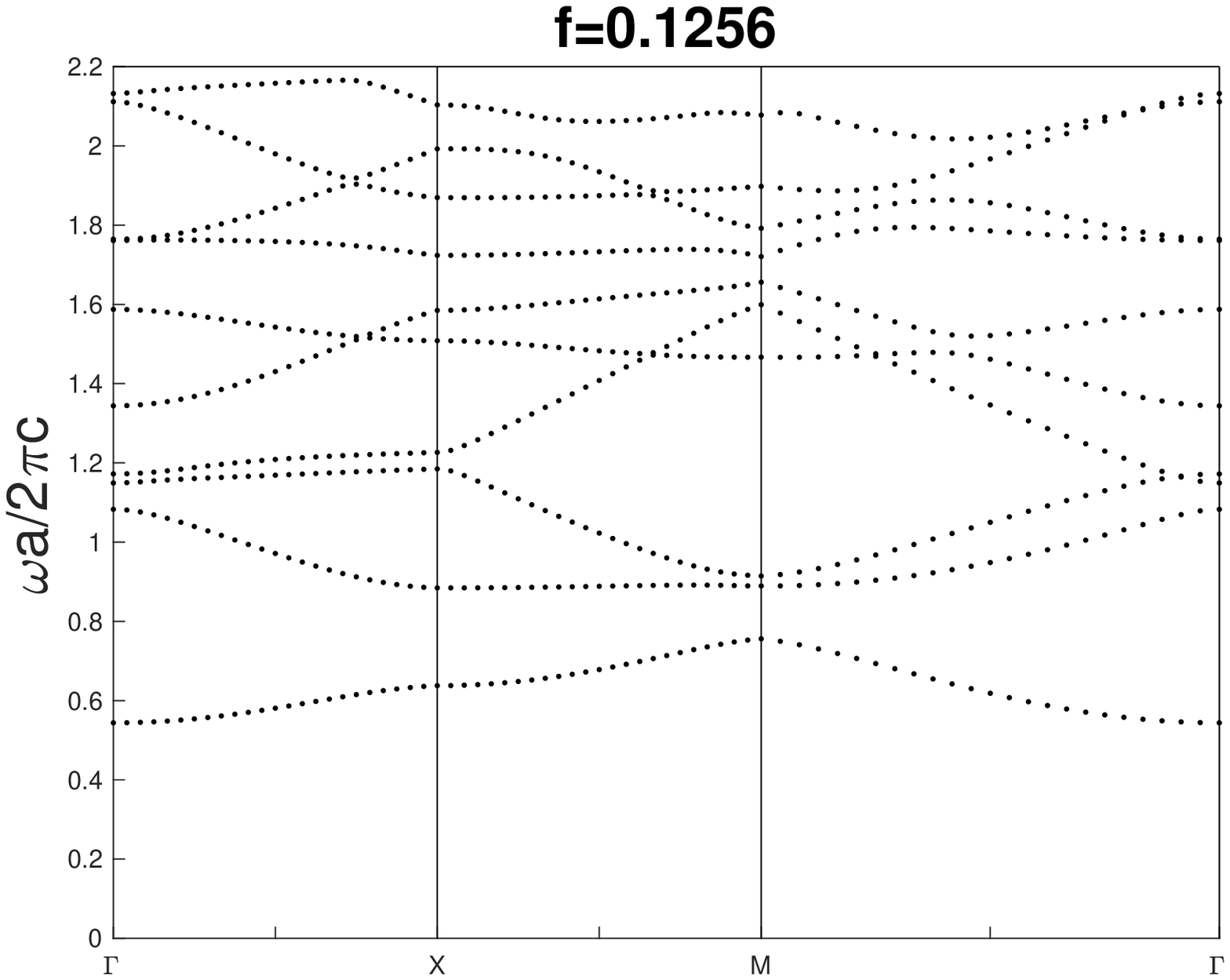}}\\
\resizebox{0.42\textwidth}{!}{\includegraphics{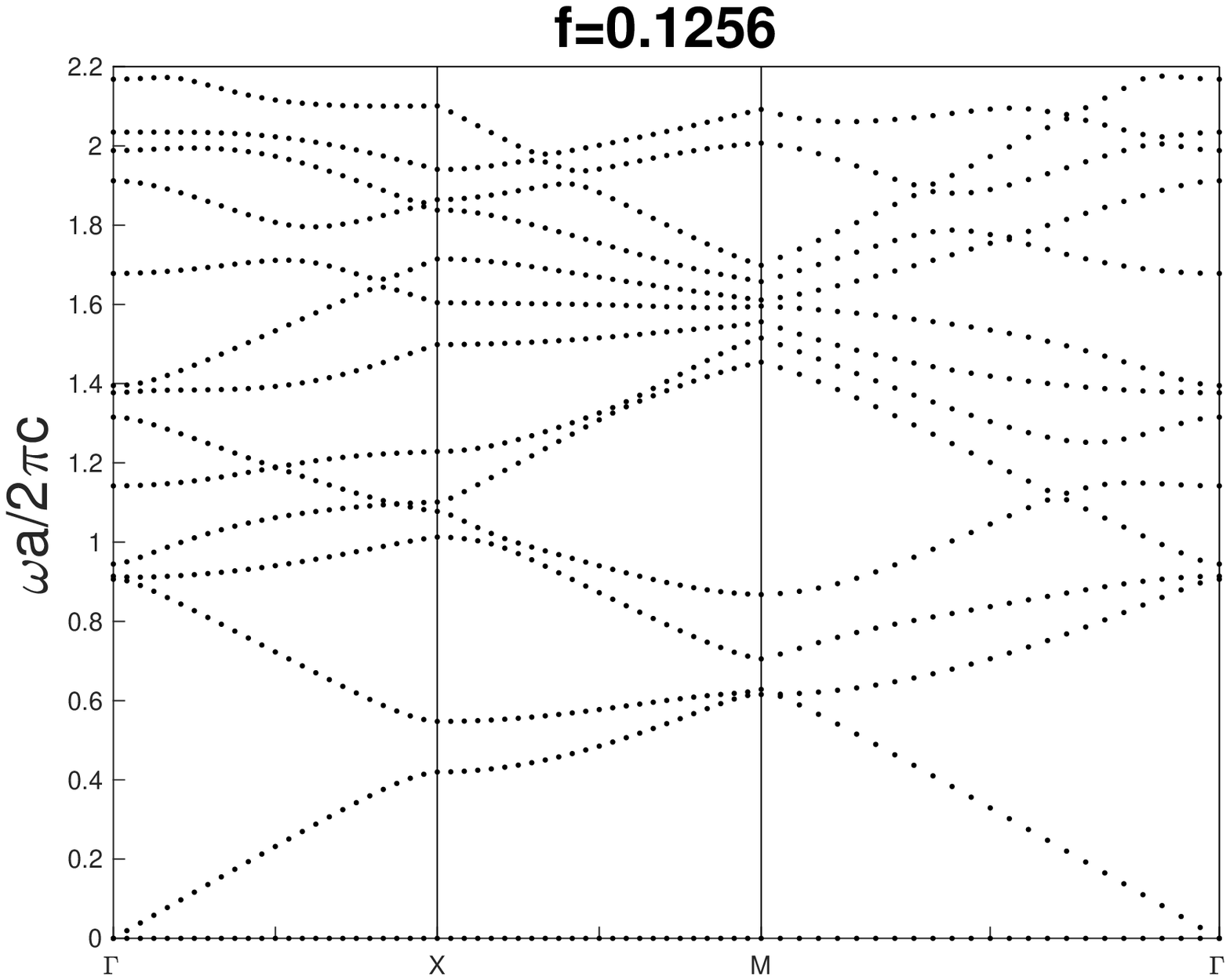}}
\end{tabular}
\end{center}
\caption{Example 1. Dispersion diagrams for metal rods with materials described by the Drude model.
Top: TE case, Bottom: TM case. }
\label{Example1}
\end{figure}

\textbf{Example 1.} The dispersive material is described by the Drude model: $\epsilon_b(\omega)=1-\frac{\omega_p^2}{\omega^2-i \omega\omega_\tau}$ where
$\omega_p = 2\pi \times1914$THz  is the plasma frequency and $\omega_\tau = 2\pi \times8.34$THz is the damping frequency \cite{raman2010, degirmenci2013}. 
The filling fraction $f=0.1256$ such that the radius of the disc is $r \approx 0.2$.
We show the dispersion diagrams in Fig.~\ref{Example1}, which are consistent with Fig.~5 and Fig.~6 of \cite{degirmenci2013}.

\begin{figure}[h!]
\begin{center}
\begin{tabular}{lll}
\resizebox{0.42\textwidth}{!}{\includegraphics{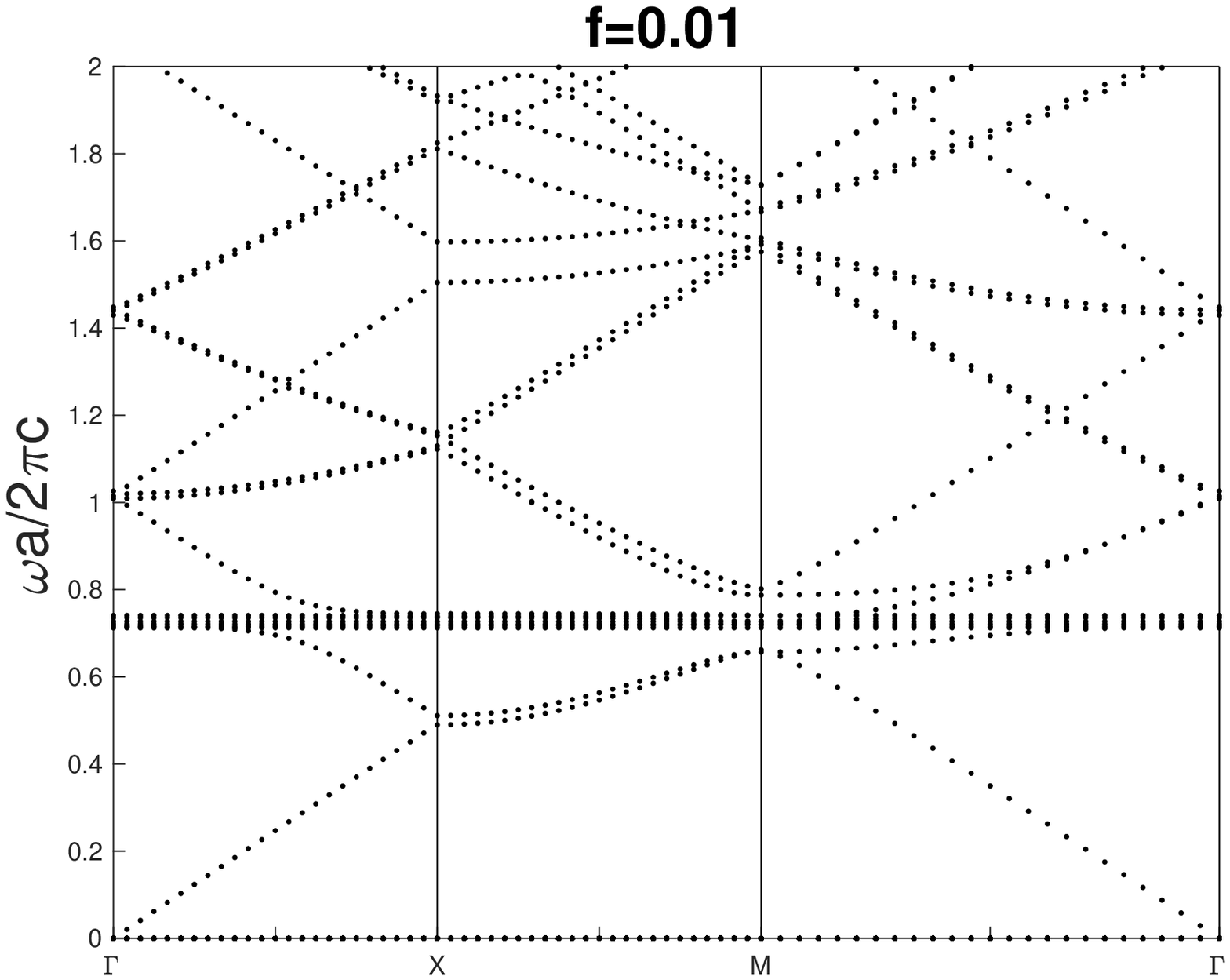}}\\
\resizebox{0.42\textwidth}{!}{\includegraphics{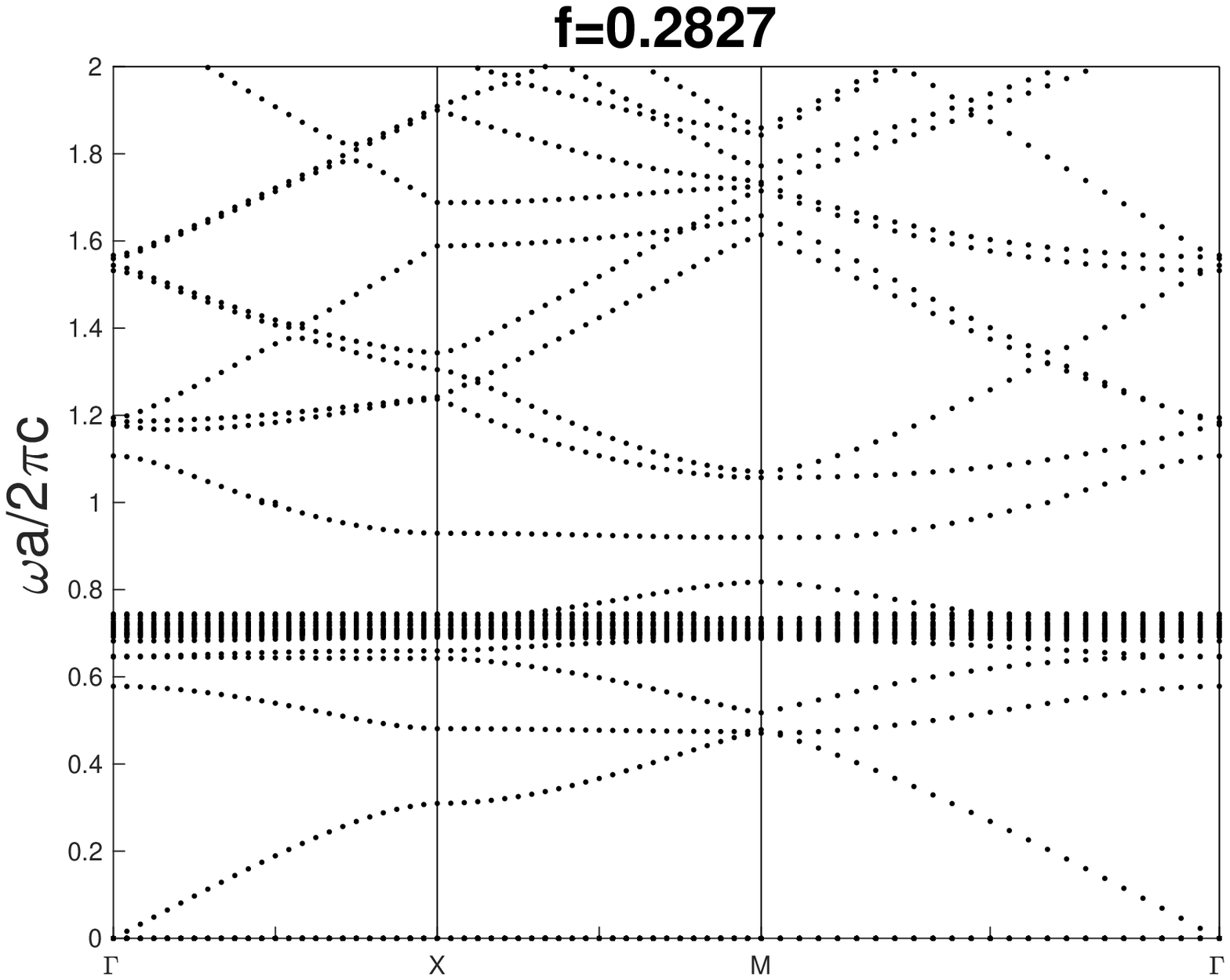}}
\end{tabular}
\end{center}
\caption{Example 2. Dispersion diagrams for metal rods of a lossy materials for the TM case.
Top: $f=0.01$, Bottom $f=0.2827$.}
\label{Example2}
\end{figure}

\textbf{Example 2.} The lossy dispersive material is described by
$\epsilon(\omega)=1-\frac{\omega_p^2}{\omega(\omega+i\gamma)}$ with
$\frac{\omega_p a}{2\pi c}=1, \gamma=0.01\omega_p$, where $\gamma$ is an inverse electronic relaxation time \cite{ito2001, kuzmiak1998}. 
In Fig.~\ref{Example2}, we show the dispersion diagrams with the filling fraction $f=0.01$ and $f=0.2827$, respectively, for the TM case.
The results for $f=0.2827$ is consistent with Fig.~5 in \cite{ito2001}. Note that \cite{ito2001} only shows a part of the diagram (from $\Gamma$ to $X$). 

\begin{thebibliography}{10}
\bibitem{Soukoulis} C.M. Soukoulis (ed.), Photonic Band Gap Materials, Springer, Dordrecht, 1996.
\bibitem{Sakoda2001} K. Sakoda, Optical Properties of Photonic Crystals. Springer, Berlin, 2001.
\bibitem{YangEtal} B. Yang, et al.,
	{\em Ideal Weyl points and helicoid surface states in artificial photonic crystal structures.}
	Science, 359 (2018), Issue 6379, 1013-1016.
\bibitem{ho1990} K.M. Ho, C.T. Chan, and C.M. Soukoulis, 
 	{\em Existence of a photonic gap in periodic dielectric structures.}
	Phys. Rev. Lett. 65 (1990), 3152-3155.	
\bibitem{qiu2000} M. Qiu and S. He, 
	{\em A nonorthogonal finite-difference time-domain method for computing the band structure of a two-dimensional photonic crystal with dielectric and metallic inclusions.}
	J. Appl. Phys. 87 (2000), 8268-8275.
\bibitem{axmann1999} W. Axmann and P. Kuchment, 
	{\em An efficient finite element method for computing spectra of Photonic and Acoustic band-gap materials I. Scalar Case.}
	 J. Comput. Phys. 150 (1999), 468-481.
\bibitem{dobson1999} D.C. Dobson,  
	{\em An efficient method for band structure calculations in 2D photonic crystals.}
	J. Comput. Phys. 149 (1999), 363-376.
\bibitem{spence2005} A. Spence and C. Poulton,  
	{\em Photonic band structure calculations using nonlinear eigenvalue techniques.}
	J. Comput. Phys. 204 (2005), 65-81.
\bibitem{kuzmiak1998} V. Kuzmiak, A. A. Maradudin, 
	Distribution of electromagnetic field and group velocities in two-dimensional periodic systems with dissipative metallic components, Physical Review B. (1998), 7230-7251.
 \bibitem{ito2001} T. Ito, K. Sakoda. Photonic bands of metallic systems. II. Features of surface plasmon polaritons. Physical Review B, 64, 045117 (2001).
\bibitem{toader2004} O. Toader and S. John, 
	{\em Photonic band gap enhancement in frequency-dependent dielectrics.}
	Physical Review E 70 (2004), 046605.
\bibitem{karma1996II} O. Karma, 
	{\em Approximation in eigenvalue problems for holomorphic Fredholm operator functions. II. (Convergence rate).} 
	Numer. Funct. Anal. Optim. 17 (1996), no. 3-4, 389-408.
\bibitem{gohberg2009} I. Gohberg and J. Leiterer, Holomorphic operator functions of one variable and applications. Birkh\"{a}user Verlag, Basel, 2009.
\bibitem{engstrom2010} C. Engstr\"{o}m, 
	{\em On the spectrum of a holomorphic operator-valued function with applications to absorptive photonic crystals.}
	Math. Models Methods Appl. Sci. 20 (2010), 1319-1341.
\bibitem{xiao2020} W. Xiao, B. Gong, J. Sun and Z. Zhang, 
	{\em A new finite element approach for the Dirichlet eigenvalue problem.}
	Appl. Math. Lett. 105 (2020), 106295.
\bibitem{huang2016} R. Huang, A. Struthers, J. Sun and R. Zhang, 
	{\em Recursive integral method for transmission eigenvalues.}
	J. Comput. Phys. 327 (2016), 830-840.
\bibitem{huang2018} R. Huang, J. Sun and C. Yang, 
	{\em Recursive integral method with Cayley transformation.}
	Numer. Linear Algebra Appl. 25 (2018), no. 6, e2199.
\bibitem{Huang2020} R. Huang, J. Sun and C. Yang,
	{\em A multilevel spectral indicator method for eigenvalues of large non-Hermitian matrices.}
	CSIAM Trans. Appl. Math., accepted, 2020. arXiv:2006.16117.
 \bibitem{SunZhou2016} J. Sun and A. Zhou, Finite element methods for eigenvalue problems. 
	CRC Press, Taylor \& Francis Group, Boca Raton, 2016.
\bibitem{raman2010} A. Raman and S. Fan, 
	{\em Photonic band structure of dispersive metamaterials formulated as a Hermitian eigenvalue problem.}
	Phys. Rev. Lett. 104 (2010), 087401.




\bibitem{kuchment1993} P. Kuchment,  Floquet Theory for Partial Differential Equations. Birkh\"{a}user Verlag, Basel (1993).
\bibitem{degirmenci2013}
E. Degirmenci, P. Landais, Finite element method analysis of band gap and transmission of two-dimensional metallic photonic crystals at terahertz frequencies, 
Applied Optics. 52, 7367-7375 (2013).
\end{thebibliography}
\end{document}